\documentclass[12pt]{amsart}

\usepackage{vmargin}
\usepackage[dvips]{graphicx}
\usepackage{color}
\input{amssym.def}
\input{amssym}
\usepackage{a4wide}
\usepackage{supertabular}
\usepackage{array}
\usepackage{multicol}

\setmarginsrb{3cm}{2cm}{3cm}{3cm}{75pt}{20pt}{20pt}{30mm}
\def\derpar#1#2{\frac{\partial#1}{\partial#2}}
\def\R{\mathbb R}
\def\N{\mathbb N}
\def\T{\mathbb T}
\def\C{\mathbb C}

\def\Z{\mathbb Z}

\def\var{\varepsilon}

\def\pa{\partial}

\def\Om{\Omega}

\def\ov{\overline}
\def\cal{\mathcal}

\def\hat{\widehat}

\def\tilde{\widetilde}

\def\Id{{\rm Id}\,}
\def\dps{\displaystyle}

\def\<{\langle}
\def\>{\rangle}

\def\cF{{\cal F}} 

\def\cZ{{\cal Z}}

\def\med{\medskip}
\def\sm{\smallskip}
\def\bul{$\bullet$\ }

\def\begeq{\begin{equation}}
\def\endeq{\end{equation}}
\def\begar{\begin{eqnarray}}
\def\endar{\end{eqnarray}}
\def\begar*{\begin{eqnarray*}}
\def\endar*{\end{eqnarray*}}
\def\begal{\begin{align}}
\def\endal{\end{align}}
\def\begal*{\begin{align*}}
\def\endal*{\end{align*}}

\newtheorem{Thm}{Theorem}

\theoremstyle{definition}

\theoremstyle{remark}

\newtheorem*{Thm*}{Theorem}
\newtheorem*{Lem*}{Lemma}
\newtheorem*{Conj*}{Conjecture}
\newtheorem*{Cor*}{Corollary}
\newtheorem*{Def*}{Definition}
\newtheorem*{Prop*}{Proposition}
\newtheorem*{Exo*}{Exercise}
\newtheorem*{Exs*}{Examples}
\newtheorem*{Ex*}{Example}
\newtheorem*{Rk*}{Remark}
\newtheorem*{Rks*}{Remarks}

\def\signcv{\bigskip \begin{center} {\sc C\'edric Villani\par\vspace{3mm}
ENS Lyon \& Institut Universitaire de France\par
UMPA, UMR CNRS 5669\par
46 all\'ee d'Italie\par
69364 Lyon Cedex 07\par
FRANCE\par\vspace{3mm}
e-mail:} \tt{cvillani@umpa.ens-lyon.fr} \end{center}}

\def\signcm{\bigskip \begin{center} {\sc Cl\'ement
      Mouhot\par\vspace{3mm}
University of Cambridge\par
DAMTP, Centre for Mathematical Sciences\par
Wilberforce Road\par
Cambridge CB3 0WA\par
ENGLAND\par
{\it On leave from:}\par
\'ENS Paris \& CNRS\par
DMA, UMR CNRS 8553\par
45 rue d'Ulm\par
F 75320 Paris cedex 05\par
FRANCE\par\vspace{3mm}
e-mail:} \tt{Clement.Mouhot@ens.fr} \end{center}}

\begin{document}

\title{Landau damping}

\author{C. Mouhot}
\author{C. Villani}

\maketitle

\vspace*{-10mm}

\begin{abstract} 
  In this note we present the main results from the recent work
  \cite{MV}, which for the first time establish Landau damping in a
  nonlinear context.
\end{abstract}



\section{Introduction}
\label{sec:intro}

The ``standard model'' of classical plasma physics is the
Vlasov--Poisson--Landau equation \cite{vlasov,landau}, here written
with periodic boundary conditions and in adimensional units:
\begeq\label{VL} \derpar{f}{t} + v\cdot\nabla_x f + F[f]\cdot\nabla_v
f = \frac{\log\Lambda}{2\pi \Lambda}\, Q_L(f,f),
\endeq 
where $f=f(t,x,v)$ is the electron distribution function
($t\geq 0$, $v\in\R^3$, $x\in\T^3 = \R^3/\Z^3$), \begeq\label{Ff}
F[f](t,x) = -\iint \nabla W(x-y)\, f(t,y,w)\,dw\,dy
\endeq
is the self-induced force, $W(x) = 1/|x|$ is the Coulomb interaction
potential, and $Q_L$ is the Landau collision operator, described for
instance in \cite{vill:handbook:02}.  The parameter $\Lambda$ is very
large, ranging typically from $10^2$ to $10^{30}$.

On very large time scales (say $O(\Lambda/\log\Lambda)$), dissipative
phenomena play a nonnegligeable role, and the entropy increase is
supposed to force the (slow) convergence to a Maxwellian distribution.
Thanks to the recent progress on hypocoercivity, this mechanism is now
rather well understood, as soon as global smoothness estimates are
available (see \cite{vill:hypoco} and the references therein).

Ten years after devising this collisional scenario, Landau
\cite{landau} formulated a much more subtle prediction: the stability
of homogeneous equilibria satisfying certain conditions --- for
instance any function of $|v|$, not necessarily Gaussian --- on much
shorter time scales (say $O(1)$), by means of purely conservative
mechanisms.  This phenomenon, called {\bf Landau damping}, is a
property of the (collisionless) Vlasov equation, obtained by setting
$\Lambda=\infty$ in \eqref{VL}. This is a theoretical cornerstone of
the classical plasma physics (among a large number of references let
us mention \cite{akhiezer}). Similar damping phenomena also occur in
other domains of physics.

The Landau damping has been since long understood at the linearized
level \cite{degond:landau,MF:landaudamping,saenz}, but the study of
the full (nonlinear) equation poses important conceptual and technical
problems. As a consequence, up to now the only existing results were
proving existence of {\em some} damped solutions with prescribed behavior
as $t\to\pm\infty$ \cite{caglimaff:VP:98,HV:landau}. 
We fill this gap in a recent work \cite{MV}, whose main
result we shall now describe.

\section{Main result}
\label{sec:resultat}

If $f$ is a function defined on $\T^d\times\R^d$, we note, for any
$k\in\Z^d$ and $\eta\in\R^d$,
\[ \hat{f}(k,v) = \int_{\T^d} f(x,v)\, e^{-2i\pi k\cdot x}\,dx,\qquad
\tilde{f}(k,\eta) = \iint_{\T^d \times \R^d} 
   f(x,v)\,e^{-2i\pi k\cdot x}\,e^{-2i\pi \eta\cdot v}\,dv\,dx.\]
We also set, for $\lambda,\mu, \beta>0$,
\begin{equation}\label{norm1}
 \|f\|_{\lambda,\mu,\beta} = 
\sup_{k,\eta} 
\Bigl( |\tilde{f}(k,\eta)|\,e^{2\pi\lambda|\eta|}\, e^{2\pi \mu |k|}\Bigr)
+ \iint_{\T^d \times \R^d} |f(x,v)|\,e^{2\pi \beta|v|}\,dv\,dx.
\end{equation}

\begin{Thm}[nonlinear Landau damping for general interaction]
\label{thmlim}
Let $d\geq 1$, and $f^0:\R^d\to\R_+$ an analytic velocity profile. Let
$W:\T^d\to\R$ be an interaction potential.  For any $k\in\Z^d$,
$\xi\in\C$, we set
\[ {\cal L}(k,\xi) = -4\pi^2\, \hat{W}(k) \int_0^\infty e^{2\pi
  |k|\xi^*t}\,|\tilde{f}^0(kt)|\,|k|^2\,t\,dt.\] We assume that there
is $\lambda>0$ such that, for $\epsilon$ small enough,
\begeq\label{condf0} \sup_{\eta\in\R^d}\ |\tilde{f}^0(\eta)|\,
e^{2\pi\lambda|\eta|} \leq C_0,\qquad \sum_{n\in\N^d}
\frac{\lambda^n}{n!} \|\nabla_v^n f^0\|_{L^1(dv)} \leq C_0,
\endeq
\begeq\label{condL} \inf_{k\in\Z^d}\ \inf_{0 \leq \, \Re \, \xi \,
  < \lambda}\ \bigl|{\cal L}(k,\xi) - 1 \bigr| \geq\kappa>0
\endeq
\begeq\label{condW} \exists \, \gamma\geq 1;\ \forall \, k\in
\Z^d;\qquad |\hat{W}(k)| \leq \frac{C_W}{|k|^{1+\gamma}}.
\endeq
Then as soon as $0<\lambda'<\lambda$, $0<\mu'<\mu$, $\beta>0$,
$r\in\N$, there are $\var>0$ and $C>0$, depending on
$d,\gamma,\lambda,\lambda',\mu,\mu',C_0,\kappa,C_W,\beta,r$, such that
if $f_i\geq 0$ satisfies \begeq\label{condfi} \delta:=
\|f_i-f^0\|_{\lambda,\mu,\beta} \leq \var,
\endeq
then the unique solution of the nonlinear Vlasov equation
\begeq\label{nlv}
\derpar{f}{t} + v\cdot\nabla_x f + F[f]\cdot\nabla_v f = 0,\qquad
F[f](t,x) = -\iint \nabla W(x-y) \, f(t,y,w)\,dw\,dy,
\endeq
defined for all times and such that $f(0,\,\cdot\,)=f_i$, satisfies
\begeq\label{rr}
\bigl\|\rho(t,\,\cdot\,)-\rho_\infty\bigr\|_{C^r(\T^d)} \leq
C\,\delta\,e^{-2\pi \lambda'|t|},
\endeq
where $\rho(t,x) = \int f(t,x,v)\,dv$, $\rho_\infty = \iint
f_i(x,v)\,dv\,dx$. Futhermore, there are analytic profiles
$f_{+\infty}(v)$, $f_{-\infty}(v)$ such that
\[ f(t,\,\cdot\,) \xrightarrow[]{t\to\pm\infty} f_{\pm\infty} \qquad
\text{weakly}\]
\[ \int f(t,x,\,\cdot\,)\,dx \xrightarrow[]{t\to\pm\infty}
f_{\pm\infty}\qquad \text{strongly (in $C^r(\R^d_v)$)},\]
these convergences being also $O(\delta\,e^{-2\pi\lambda'|t|})$.
\end{Thm}
\med

This theorem, entirely constructive, is almost optimal, as the
following comments show. \med

\noindent{\bf Comments on the assumptions:} The periodic boundary
conditions of course are debatable; in any case, the counterexamples
of Glassey and Schaeffer \cite{GS} show that some confinement
mechanism --- or at least a limitation on the wavelength --- is
mandatory. Condition \eqref{condf0} quantitatively expresses the
analyticity of the profile $f^0$, without which we could not hope for
an exponential convergence. The inequality \eqref{condL} is a linear
stability condition, roughly optimal, covering all physically
interesting cases: in particular the (attractive) Newton interaction
for wavelengths shorter than the Jeans unstability length; and the
(repulsive) Coulomb interaction around radially symmetric profiles
$f^0$, for all wavelengths. On the other hand, condition \eqref{condW}
shows up only in the nonlinear stability; it is satisfied by
Coulomb and Newton interactions as a limit case. As for the condition
\eqref{condfi}, its perturbative nature is natural in view of
theoretical speculations and numerical studies in the subject. \med

\med

\noindent{\bf Comments on the conclusions:}

\begin{enumerate}
\item The large-time convergence is based on a reversible, purely
  deterministic mechanism, without any Lyapunov functional neither
  variational interpretation. The asymptotic profiles $f_{\pm\infty}$
  eventually keep the memory of the initial datum and the
  interaction. This convergence ``for no reason'' was not really
  expected, since the quasilinear theory of Landau damping
  \cite[Vol. II, Section 9.1.2]{akhiezer} predicts convergence only
  after taking average on statistical ensembles.

\item This result can be interpreted in the spirit of the KAM theorem:
  for the linear Vlasov equation, convergence is forced by an infinity
  of invariant subspaces, which make the model ``completely integrable'';
  as soon as one adds a nonlinear coupling, the invariance goes away
  but the convergence remains. 

\item Given a stable equilibrium profile $f^0$, we see that an entire
  neighborhood --- in analytic topology --- of $f^0$ is filled by
  homoclinic or (in general) heteroclinic trajectories. Only infinite
  dimension allows this remarkable behavior of the nonlinear Vlasov
  equation.

\item The large time convergence of the distribution function holds
  only in the weak sense; the norms of velocity derivatives grow
  quickly in large time, which reflects a filamentation in phase
  space, and a transfer of energy (or information) from low to high
  frequencies (``weak turbulence'').

\item It is this transfer of information to small scales which allows
  to reconcile the reversibility of the Vlasov--Poisson equation with
  the seemingly irreversible large-time behavior. Let us note that the
  ``dual'' mechanism of transfer of energy to large scales, also
  called {\em radiation}, was extensively studied in the setting of
  Hamiltonian systems.
\end{enumerate}

Much more comments, both from the mathematical and the physical sides,
can be found in \cite{MV}.

\section{Linear stability}
\label{sec:lin}

The linear stability is the first step of our study; it only requires
a reduced technical investment.

After linearization around a homogeneous equilibrium $f^0$, the Vlasov
equation becomes \begeq\label{LV} \derpar{h}{t} + v\cdot\nabla_x h -
(\nabla W\ast \rho)\cdot\nabla_v f^0 =0,\qquad \rho = \int h\,dv.
\endeq
It is well-known that this equation decouples into an infinite number
of independent equations governing the modes of $\rho$: for all
$k\in\Z^d$ and $t\geq 0$, \begeq\label{rho} \hat{\rho}(t,k) - \int_0^t
K^0(t-\tau,k)\,\hat{\rho}(\tau,k)\,d\tau = \tilde{h}_i(k,kt),
\endeq
where $h_i$ is the initial datum, and $K^0$ an integral kernel
depending on $f^0$: \begeq\label{K0} K^0(t,k) = -
4\pi^2\,\hat{W}(k)\,\tilde{f}^0(kt)\,|k|^2\,t.
\endeq

Then from classical results on Volterra equations we deduce that for
all $k\neq 0$ the decay of $\hat{\rho}(t,k)$ as $t\to\infty$ is
essentially controlled by the worst of two convergence rates: \sm

\bul the convergence rate of the source term in the right-hand side of
\eqref{rho}, which depends only on the regularity of the initial datum
in the velocity variable; \sm

\bul $e^{-\lambda t}$, where $\lambda$ is the largest positive real
number such that the Fourier--Laplace transform (in the $t$ variable)
of $K^0$ does not approach the value~1 in the strip $\{0\leq \Re z
\leq \lambda\}\subset \C$. The problem lies in finding sufficient
conditions on $f^0$ to guarantee the strict positivity of $\lambda$.
\sm

Since Landau, this study is traditionally performed thanks to the
Laplace transform inversion formula; however, with a view to the
nonlinear study, we prefer a more elementary and constructive
approach, based on the plain Fourier inversion formula.

With this method we establish the linear Landau damping, under
conditions \eqref{condL} and \eqref{condf0}, for any interaction $W$
such that $\nabla W \in L^1(\T^d)$, and any analytical initial
condition (without any size restriction in this linear context). We
recover as particular cases all the results previously established on
the linear Landau damping \cite{degond:landau,MF:landaudamping,saenz};
but we also cover for instance Newton interaction. Indeed, condition
\eqref{condL} is satisfied as soon as {\em any one} of the following
conditions is satisfied: \sm

(a) $\forall \, k\in\Z^d$, $\forall \, z\in\R$, $\hat{W}(k)\geq 0$, \
$z\phi'_k(z)\leq 0$, where $\phi_k$ is the ``marginal'' of $f^0$ along
the direction $k$, defined by
\[ \phi_k(z) = \int_{\frac{kz}{|k|}+k^\bot} f^0(w)\,dw;\]

(b) $\dps 4\pi^2\, \left(\max_{k \not = 0} \,|\hat{W}(k)|\right)
\left(\sup_{|\sigma|=1} \int_0^\infty |\tilde{f}^0(r\sigma)|\,r\,dr
\right) < 1$.  \med

We refer to \cite[Section 3]{MV} for more details.

\section{Nonlinear stability}
\label{sec:nonlin}

To establish the nonlinear stability, we start by introducing analytic
norms which are ``hybrid'' (based on the size of derivatives in the
velocity variable, and on the size of Fourier coefficients in the
position variable) and ``gliding'' (the norm will change with time to
take into account the transfer to small velocity scales). Five indices
provide all the necessary flexibility: \begeq\label{Z}
\|f\|_{\cZ^{\lambda,(\mu,\gamma);p}_\tau} = \sum_{k\in\Z^d}
\sum_{n\in\N^d} e^{2\pi\mu|k|}\, (1+|k|)^\gamma\,
\frac{\lambda^n}{n!}\, \Bigl\|\bigl(\nabla_v + 2i\pi \tau k\bigr)^n
\hat{f}(k,v)\Bigr\|_{L^p(dv)}.
\endeq
(By default $\gamma=0$.) A tedious injection theorem ``\`a la
Sobolev'' compares these norms to more traditional ones, such as the
$\|f\|_{\lambda,\mu,\beta}$ norms appearing in \eqref{norm1}.

The $\cZ$ norms enjoy remarkable properties with respect to
composition and product.  The parameter $\tau$ partly compensates for
filamentation.  Finally, the hybrid nature of these norms is well
adapted to the geometry of the problem.  If $f$ depends only on $x$,
the norm \eqref{Z} coincides with the norm
$\cF^{\lambda\tau+\mu,\gamma}$ defined by \begeq\label{F}
\|f\|_{\cF^{\lambda\tau+\mu,\gamma}} = \sum_{k\in\Z^d}
|\hat{f}(k)|\,e^{2\pi(\lambda\tau+\mu)|k|}\, (1+|k|)^\gamma.
\endeq
(We also use the ``homogeneous'' version
$\dot{\cF}^{\lambda\tau+\mu,\gamma}$ where the mode $k=0$ is removed.)

Then the Vlasov equation is solved by a Newton scheme, whose first
step is the solution of the linearized equation around $f^0$:
\[ f^n = f^0 + h^1 + \ldots + h^n, \]
\begeq\label{N1}
\begin{cases} 
\pa_t h^1 + v\cdot\nabla_x h^1 + F[h^1]\cdot\nabla_v f^0 =0 \\[2mm]
  h^1(0,\,\cdot\,) = f_i -f^0
\end{cases} \endeq
\begeq\label{N2}
n \ge 1, \quad 
\begin{cases} \pa_t h^{n+1} + v\cdot\nabla_x h^{n+1} +
  F[f^{n}]\cdot\nabla_v h^{n+1} + F[h^{n+1}]\cdot\nabla_v f^{n} =
  - F[h^{n}]\cdot\nabla_v h^{n}\\[2mm]
  h^{n+1}(0,\,\cdot\,) = 0.
\end{cases} \endeq

In a first step, we establish the short-time analytic regularity of
$h^n(\tau,\,\cdot\,)$ in the norm $\cZ^{\lambda,(\mu,\gamma);1}_\tau$;
this step, in the spirit of a Cauchy--Kowalevskaya theorem, is
performed thanks to the identity \begeq\label{leftddtnorm}
\left.\frac{d}{dt}^+ \right|_{t=\tau}
\|f\|_{\cZ^{\lambda(t),\mu(t);p}_\tau} \leq -\frac{K}{1+\tau}\,
\|\nabla f\|_{\cZ^{\lambda(\tau),\mu(\tau);p}_\tau},
\endeq
where $\lambda(t) = \lambda -Kt$, $\mu(t) = \mu - Kt$.
\sm

In a second step, we establish uniform in time estimates on each
$h^n$, now with a partly Eulerian and partly Lagrangian method,
integrating the equation along the characteristics
$(X^n_{\tau,t},V^n_{\tau,t})$ created by the force $F[f^n]$.  (Here
$\tau$ is the initial time, $t$ the current time, $(x,v)$ the initial
conditions, $(X^n,V^n)$ the current conditions.) The smoothness of
these characteristics is expressed by controls in hybrid norm on the
operators $\Om^n_{t,\tau}(x,v) =
(X^n_{t,\tau},V^n_{t,\tau})(x+v(t-\tau),v)$, which compare the
perturbed dynamics to the unperturbed one; these are informally called
scattering operators.

Then we propagate a number of estimates along the scheme; the most
important are (slightly simplifying) \begeq\label{suphn}
\sup_{\tau\geq 0} \left\| \int_{\R^d} h^n
  \bigl(\tau,\,\cdot\,,v\bigr)\,dv \right\|_{\cF^{\lambda_n\tau +
    \mu_n}}\leq \delta_n,
\endeq
\begeq\label{supttauhn} \sup_{t\geq\tau\geq 0} \Bigl\|h^n
\bigl(\tau,\Om^n_{t,\tau}\bigr)
\Bigr\|_{\cZ^{\lambda_n(1+b),\mu_n;1}_{\tau-\frac{bt}{1+b}}} \leq
\delta_n,\qquad b= b(t) = \frac{B}{1+t},
\endeq
\begeq\label{Omn} \Bigl\|\Om^n_{t,\tau}-\Id
\Bigr\|_{\cZ^{\lambda_n(1+b),(\mu_n,\gamma);\infty}_{\tau-\frac{bt}{1+b}}}
\leq C\, \left(\sum_{k=1}^n
  \frac{\delta_k\,e^{-2\pi(\lambda_k-\lambda_{n+1}) \tau}} {2\pi
    (\lambda_k-\lambda_{n+1})^2} \right)\,\min\{ t-\tau \, ; \, 1\}.
\endeq 

Notice, in \eqref{suphn}, the linear increase in the regularity of the
spatial density, which comes at the same time as the deterioration of
regularity in the $v$ variable. In \eqref{supttauhn}, the additional
time-shift in the indices by the function $b(t)$ will be crucial to
absorb error terms coming from the composition; the constant $B$
itself is determined by the previous small-time estimates. Finally, in
\eqref{Omn}, notice the uniform in $t$ control, and the improved
estimates in the limit cases $t\to\tau$ and $\tau\to\infty$; also this
is important for handling error terms. The constants $\lambda_n$ and
$\mu_n$ decrease at each stage of the scheme, converging --- not too
fast --- to positive limits $\lambda_\infty$, $\mu_\infty$; at the
same time, the constants $\delta_n$ converge extremely fast to~0,
which guarantees ``by retroaction'' the uniformity of the constants in
the right-hand side of \eqref{Omn}.

The estimates \eqref{Omn} are obtained by repeated application of
fixed point theorems in analytic norms.  Another crucial ingredient to
go from stage $n$ to stage $n+1$ is the mechanism of {\bf regularity
  extortion}, which we shall now describe in a simplified
version. Given two distribution functions $f$ and $\ov{f}$, depending
on $t,x,v$, let us define
\[ \sigma(t,x) = \int_0^t \int_{\R^d} \bigl( F[f]\cdot\nabla_v
\ov{f}\bigr) \bigl(\tau,x-v(t-\tau),v\bigr)\,dv\,d\tau.\] This
quantity can be interpreted as follows: if particles distributed
according to $f$ exert a force on particles distributed according to
$\ov{f}$, then $\sigma$ is the variation of density $\int f\,dv$
caused by the {\em reaction} of $\ov{f}$ on $f$. We show that if
$\ov{f}$ has a high gliding regularity, then the regularity of
$\sigma$ in large time is better than what would be expected:
\begeq\label{sigmaleq} \|\sigma(t,\,\cdot\,)\|_{\dot{\cF}^{\lambda
    t+\mu}} \leq \int_0^t K(t,\tau)\,
\bigl\|F\bigl[f(\tau,\,\cdot\,)\bigr]\bigr\|_{\cF^{\lambda\tau+\mu,\gamma}}\,d\tau,
\endeq
where
\[ K(t,\tau) = \left[ \sup_{0\leq s\leq t} \left(\frac{\bigl\|\nabla_v
      \ov{f}(s,\,\cdot\,)\bigr\|_{\cZ^{\ov{\lambda},\ov{\mu};1}_s}}
    {1+s}\right) \right] \,(1+\tau)\, \sup_{k\neq 0,\ \ell\neq 0}\
\frac{e^{-2\pi(\ov{\lambda}-\lambda)|k(t-\tau)+\ell\tau|}\, e^{-2\pi
    (\ov{\mu}-\mu)|\ell|}}{1+|k-\ell|^\gamma}.\] The kernel
$K(t,\tau)$ has integral $O(t)$ as $t\to\infty$, which would let us
fear a violent unstability; but it is also more and more concentrated
on discrete times $\tau=kt/(k-\ell)$; this is the effect of {\bf
  plasma echoes}, discovered and experimentally observed in the
sixties \cite{echo:expe}.  The stabilizing role of the echo
phenomenon, related to the Landau damping, is uncovered in our study.

Then we analyze the nonlinear response due to echoes.  If $\gamma>1$,
from \eqref{sigmaleq} one deduces that the response is subexponential,
and therefore can be controlled by an arbitrarily small loss of
gliding regularity, at the price of a gigantic constant, which later
will be absorbed by the ultrafast convergence of the Newton scheme.
In the end, part of the gliding regularity of $\ov{f}$ has been
converted into a large-time decay.

When $\gamma=1$, a finer strategy is needed. To handle this case, we
work on the response mode by mode, that is, estimating the size of
$\hat{\rho}(t,k)$ for all $k$, via an infinite system of
inequalities. Then we are able to take advantage of the fact that
echoes occurring at different frequencies are asymptotically rather
well separated. For instance, in dimension~1, the dominant echo
occurring at time $t$ and frequency $k$ corresponds to $\tau =
kt/(k+1)$.

In practice, straight trajectories in \eqref{sigmaleq} must be
replaced by characteristics (this reflects the fact that $\ov{f}$ also
exerts a force on $f$), which is a source of considerable technical
difficulties. Among the tools used to overcome them, let us mention a
second mechanism of regularity extortion, acting in short time and
close in spirit to velocity-averaging lemmas; here is a simplified
version of it: \begeq\label{sigmaleq2}
\|\sigma(t,\,\cdot\,)\|_{\dot{\cF}^{\lambda t+\mu}} \leq \int_0^t
\bigl\|F\bigl[f(\tau,\,\cdot\,)\bigr]\bigr\| _{\cF^{\lambda [\tau -
    b(t-\tau)] +\mu,\gamma}} \, \bigl\|\nabla f(\tau,\,\cdot\,)\bigr\|
_{\cZ^{\lambda(1+b), (\mu,0);1} _{\tau - bt/(1+b)}} \,d\tau.
\endeq
We see in \eqref{sigmaleq2} that the regularity of $\sigma$ is better
than that of $F[f]$, with a gain that degenerates as $t\to\infty$ or
$\tau\to t$. 
\med

\bibliographystyle{acm}
\bibliography{./note-eng}


\signcm
\signcv

\end{document}